\documentclass[12pt]{article}
\usepackage{amssymb}
\usepackage{amsfonts}
\usepackage{graphicx}
\usepackage{amsmath}
\usepackage{booktabs}
\usepackage{multirow}

\newtheorem{example}{Example}

\begin{document}
\date{}
\title{Error estimates for a Gaussian rule involving Bessel functions}
\author{Eleonora Denich\thanks{Dipartimento di Matematica e Geoscienze, Universit\`{a} di Trieste, Trieste, Italy, eleonora.denich@phd.units.it}}
\maketitle

\begin{abstract}
This paper deals with the estimation of the quadrature error of a Gaussian formula for weight functions involving fractional powers, exponentials and Bessel functions of the first kind. 
For this purpose, in this work the averaged and generalized averaged Gaussian rules are employed, together with a tentative a priori approximation of the error. The numerical examples confirm the reliability of these approaches.
\end{abstract}



\section{Introduction}

In this work we consider the approximation of integrals of the type
\begin{equation} \label{integrale}
I_{\nu,\alpha,c}(f) = \int_0^{\infty} f(x) x^{\alpha} e^{-cx} J_{\nu} (x) dx,
\end{equation}
where $J_{\nu}$ is the Bessel function of the first kind of order $\nu \geq 0$ (see \cite{Bessel} for a background), $\alpha>-1$, $c>0$ and $f$ is a smooth function. 
We point out that integrals of type (\ref{integrale}) are strongly related to Hankel transforms, which commonly appear in problems of mathematical physics and applied mathematics having axial symmetry.
An example of application arises in geophysical electromagnetic survey. In particular, electromagnetic fields over a layered earth due to magnetic dipoles above the surface can be expressed in integral form as in (\ref{integrale}) (see \cite{FastEM}).

Since $\left\vert J_{\nu} (x)\right\vert \leq 1$, for $\nu \geq 0$, $x \in \mathbb{R}$ (see \cite[p.362]{Abramowitz}), a Gaussian quadrature rule for the computation of integrals of type (\ref{integrale}) was constructed in \cite{GaussBessel} by rewriting (\ref{integrale}) as  
\begin{equation*}
I_{\nu,\alpha,c}^J(f)-I_{\alpha,c}^L(f),
\end{equation*}
with
\begin{equation} \label{integrale1}
I_{\nu,\alpha,c}^J(f):=\int_0^{\infty} f(x) x^{\alpha} e^{-cx} [J_{\nu} (x)+1] dx,
\end{equation}
and	
\begin{equation} \label{integrale2}
I_{\alpha,c}^L(f):=\int_0^{\infty} f(x) x^{\alpha} e^{-cx} dx.
\end{equation}
In this setting, the authors considered the approximations
\begin{equation*} 
I_{\nu,\alpha,c}^J(f) \approx I_n^J(f) \quad {\rm and} \quad I_{\alpha,c}^L(f) \approx I_n^L(f),
\end{equation*}
where
\begin{equation} \label{GJ}
I_n^J(f) = \sum_{i=1}^n w_i^{(n)} f\left( x_i^{(n)} \right)
\end{equation}
is the Gaussian rule relative to the weight function
\begin{equation} \label{pesoBessel}
    w_{\nu, \alpha,c}(x):=x^{\alpha} e^{-cx} [J_{\nu}(x)+1] \quad {\rm on} \quad [0, +\infty),
\end{equation}
and
\begin{equation} \label{GL}
I_n^L(f) = \sum_{i=1}^n \lambda_i^{(n)} f\left(\xi_i^{(n)} \right)
\end{equation}
is a slight modification of the Gauss-Laguerre quadrature rule, i.e., relative to the weight function 
\begin{equation} \label{pesoLagC}
    w_{\alpha,c}(x) = x^{\alpha} e^{-cx} \quad {\rm on} \quad [0, +\infty].
\end{equation}
Clearly, it is possible to consider a change of variable in (\ref{integrale2}) in order to work with the standard Laguerre rule, but our choice allows to treat in the same way integrals (\ref{integrale1}) and (\ref{integrale2}). 

Denoting by $I_n(f)$ the $n$-point Gaussian rule for the integral
\begin{equation*}
I(f)= \int_0^{+\infty} f(x) w(x) dx,
\end{equation*}
in which $w$ is a generic weight function, it is not easy, in general, to derive an accurate estimate of the error
\begin{equation} \label{errore}
E_n(f)= I(f) - I_n(f).
\end{equation}
A classical approach is to consider the $(2n+1)$-point Gauss-Kronrod quadrature rule associated with the $n$-point Gaussian formula $I_n(f)$ (see \cite{K,G1,N}). However, in \cite{KM} the nonexistence of Gauss-Kronrod rules, for $n >2$, with real nodes and positive weights for the Gauss-Laguerre formula was proved.
As consequence, this approach is not suitable for our case.

An alternative approach was proposed by Laurie \cite{Laurie}, who introduced the so-called anti-Gaussian quadrature rule $A_{n+1}$, corresponding to $I_n$.
It is a $(n+1)$-point formula of degree $2n-1$ which integrates polynomial of degree up to $2n+1$ with an error equal in magnitude but of opposite sign to one of the $n$-point Gaussian formula.
Then the idea is to estimate error (\ref{errore}) as
\begin{equation*}
    E_n(f) \approx \tilde{A}_{2n+1}(f)-I_n(f),
\end{equation*}
where the $(2n+1)$-point formula
\begin{equation*}
\tilde{A}_{2n+1}(f) = \frac{1}{2}\left( A_{n+1}(f)+I_n(f) \right)
\end{equation*}
is commonly named averaged Gaussian formula.
The anti-Gaussian rule always exists, it is guaranteed to have real nodes and positive weights, and at most one of the nodes may be outside the integration interval (see \cite{Laurie}). 

A more general formula, given by
\begin{equation} \label{averaged_Gauss}
H_{2n+1}(f)=\frac{1}{2+\gamma} \left( (1+\gamma) I_n(f)+A_{n+1}(f) \right), \quad \gamma>0,
\end{equation}
was considered by Enrich in \cite{E}. In particular, he constructed (\ref{averaged_Gauss}), for the Laguerre and Hermite weight functions, with the parameter $\gamma$ chosen to reach the highest degree of exactness, that is, $2n+1$. We refer to this formula as the generalized averaged Gaussian rule.

More recently, Spalević \cite{S} derived a parameter free method for constructing generalized averaged Gauss rules for any weight function for which all moments exist. We denote these quadrature formulas by $\hat{A}_{2n+1}$. They have degree of exactness at least $2n+1$ and are guaranteed to have real nodes and positive weights.
Then the error is estimated as
\begin{equation} \label{E}
    E_n(f) \approx \hat{A}_{2n+1}(f)-I_n(f).
\end{equation}

In this paper the averaged Gauss rules $\tilde{A}_{2n+1}^J$, $\tilde{A}_{2n+1}^L$ and generalized averaged Gauss rules $\hat{A}_{2n+1}^J$, $\hat{A}_{2n+1}^L$, corresponding to $I_n^J$ and $I_n^L$, respectively, are constructed and used to estimate the error
\begin{equation} \label{errore totale}
E_{n, \nu,\alpha,c}(f)=E_{n, \nu,\alpha,c}^J(f)-E_{n,\alpha,c}^L(f),
\end{equation}
where
\begin{equation} \label{11bis}
    E_{n, \nu,\alpha,c}^J(f) = I_{\nu,\alpha,c}^J(f) -I_n^J(f) \quad {\rm and} \quad  E_{n,\alpha,c}^L(f) = I_{\alpha,c}^L(f) -I_n^L(f).
\end{equation}
Unfortunately, the use of the heavy notation $E_{n, \nu,\alpha,c}^{(\cdot)}$ is necessary to avoid confusion with some general results reported in the paper.

The averaged and generalized averaged Gaussian rules are easy to construct and, moreover, typically lead to quite accurate estimates of (\ref{E}) (see \cite{RS22}).
However, sometimes it can be useful to have at disposal an a priori estimate of the error to have an idea of the number of points necessary to reach a prescribed accuracy.
In this view, here we also present a tentative a priori approximation of the quadrature error.
In particular, by interpreting $w_{\nu,\alpha,c}$ as a perturbation of the weight function of the Gauss-Laguerre rule, the idea is to employ a result due to Barrett \cite{B} and relative to the asymptotic behavior of the error of the Gauss-Laguerre formula, to estimate $E_{n, \alpha,c}^L$ and $E_{n, \nu,\alpha,c}^J$.
Moreover, similarly to the Gauss-Laguerre rule, it can be verified that the weights of $I_n^J$ decay exponentially. Hence, having at disposal a reliable error estimate, a truncated approach can also be introduced, but not considered in the present paper.

Throughout this work we use the symbol $\approx$ to indicate a generic approximation.
The symbol $\sim$ is used to express the asymptotic equality.

This paper is organized as follows.
Section \ref{section2} and \ref{section3} deals with the representation of averaged and generalized averaged Gaussian rules. Moreover, some theoretical and experimental properties of the quadrature formulas are described.
In Section \ref{section4} a tentative error approximation, which allows to have an a priori estimate of the error behavior, is presented.
Section \ref{section5} deals with the discussion of some numerical examples, which show the performances of the error estimates. 
Concluding remarks can be found in Section \ref{section6}.

\section{Construction of averaged Gaussian rules} \label{section2}

Let $w$ be a generic weight function on $[0, + \infty)$ and consider the corresponding Gaussian quadrature formula
\begin{equation} \label{gaussianRule}
I_n(f) = \sum_{i=1}^n \sigma_i^{(n)} f\left( \tau_i^{(n)} \right),
\end{equation}
of degree $2n-1$, for the integral
\begin{equation*}
    I(f) = \int_0^{+ \infty} f(x) w(x) dx.
\end{equation*}
Denoting by $\alpha_k \in \mathbb{R}$, $\beta_k >0$ the recursion coefficients for the sequences of monic orthogonal polynomials $\left\lbrace \pi_j \right\rbrace_{j \geq 0}$ relative to the weight function $w$, the Gaussian rule (\ref{gaussianRule}) can be associated with the symmetric tridiagonal matrix
\begin{equation} \label{matriceJ_n}
J_n=    \begin{bmatrix}
\alpha_0 & \sqrt{\beta_1} & & & 0\\
\sqrt{\beta_1} & \alpha_1 & \sqrt{\beta_2} \\
 & \sqrt{\beta_2} & \alpha_2 & \ddots \\
 & & \ddots & \ddots & \sqrt{\beta_{n-1}}\\
0 & & & \sqrt{\beta_{n-1}} & \alpha_{n-1}
\end{bmatrix}  \in \mathbb{R}^{n \times n}.
\end{equation}
It is well known that the eigendecomposition of the above matrix provides the nodes $\tau_i^{(n)}$ and the weights $\sigma_i^{(n)}$, $i=1,\ldots,n,$ of the Gaussian rule $I_n$ (see e.g. \cite{DR} and the references therein).
Now, the corresponding $(n+1)$-point anti-Gaussian quadrature formula
\begin{equation} \label{A}
A_{n+1}(f) := \sum_{i=1}^{n+1} \tilde{\sigma}_i^{(n+1)} f\left(\tilde{\tau}_i^{(n+1)} \right),    
\end{equation}
is such that 
\begin{equation} \label{proprietà}
    I(p) - A_{n+1}(p) = - \left(I(p)-I_n(p) \right), \quad \forall p \in \mathbb{P}_{2n+1},
\end{equation}
where $\mathbb{P}_{2n+1}$ denotes the space of polynomials of degree at most $2n+1$. 
Laurie \cite{Laurie} showed that formula (\ref{A}) is associated with the symmetric tridiagonal matrix $\tilde{J}_{n+1}\in \mathbb{R}^{(n+1) \times (n+1)}$ defined by
\begin{equation} \label{Jtilde}
\tilde{J}_{n+1}=
\begin{bmatrix} 
J_n & e_n \sqrt{2 \beta_n} \\
e_n^T \sqrt{2 \beta_n} & \alpha_n  \\
\end{bmatrix}, 
\end{equation}
where $e_n=(0,\ldots,0,1)^T \in \mathbb{R}^n$.
Therefore, having at disposal the recurrence coefficients $\alpha_k$ and $\beta_k$, it is trivial to compute the nodes $\tilde{\tau}_i^{(n+1)}$ and the weights $\tilde{\sigma}_i^{(n+1)}$, $i=1,\ldots,n+1$, of the anti-Gaussian rule.
Moreover, quadrature formula (\ref{A}) has the following properties (see \cite[Theorem 1]{Laurie}):
\begin{enumerate}
    \item $\tilde{\sigma}_i^{(n+1)}>0$, $i=1,\ldots,n+1$; 
    \item the nodes $\tilde{\tau}_i^{(n+1)}$, $i=1,\ldots, n+1$, are all real and are interlaced by those of the Gaussian formula $I_n$, that is,
    \begin{equation*}
        \tilde{\tau}_1^{(n+1)}< \tau_1^{(n)}<\tilde{\tau}_2^{(n+1)}<\ldots<\tau_n^{(n)}<\tilde{\tau}_{n+1}^{(n+1)};
    \end{equation*}
    \item $\tilde{\tau}_i^{(n+1)} \in [0, +\infty)$, for $i \geq 2$; 
    \item $\tilde{\tau}_1^{(n+1)} \in [0,+\infty)$ if and only if
    \begin{equation} \label{internal}
        \frac{\pi_{n+1}(0)}{\pi_{n-1}(0)} \geq \beta_n, \quad n \geq 1.
    \end{equation} 
\end{enumerate}
At this point, the averaged quadrature formula $\tilde{A}_{2n+1}$ is defined as
\begin{equation*}
    \tilde{A}_{2n+1}(f):= \frac{1}{2} \left(I_n(f)+A_{n+1}(f)\right).
\end{equation*}
From property (\ref{proprietà}) it follows that the degree of exactness of $\tilde{A}_{2n+1}$ is $2n+1$, and the quadrature error can be estimated by
\begin{equation*}
    E_n(f) \approx \tilde{A}_{2n+1}(f)-I_n(f) = \frac{1}{2}\left(A_{n+1}(f)-I_n(f)\right).
\end{equation*}

Going back to our case, we denote by $\alpha_k^J$, $\beta_k^J$ and $\alpha_k^L$, $\beta_k^L$ the recursion coefficients for the sequences of monic orthogonal polynomials relative to $w_{\nu, \alpha,c}$ and $w_{\alpha,c}$ (see (\ref{pesoBessel})-(\ref{pesoLagC})), respectively, and by $\tilde{J}_{n+1}^J$, $\tilde{J}_{n+1}^L$ the associated tridiagonal matrices of type (\ref{Jtilde}).
In particular, we have that $\alpha_k^L$, $\beta_k^L$ are strongly related to the recurrence coefficients $a_k$, $b_k$ of the standard generalized Gauss-Laguerre rule. Indeed, it can be easily verified that
\begin{equation} \label{ab}
    \alpha_k^L = \frac{a_k}{c} \quad {\rm and} \quad \beta_k^L = \frac{b_k}{c^2}.
\end{equation}
Moreover, the monic polynomials $\left\lbrace L_k^{(\alpha,c)} \right\rbrace_{k \geq 0}$ defined by
\begin{equation} \label{polinomi Lag c}
    L_k^{(\alpha,c)}(x) = \frac{1}{c^k} {L}_k^{(\alpha)}(cx),
\end{equation}
where ${L}_k^{(\alpha)}(t)$ is the monic generalized Laguerre polynomial of degree $k$, are orthogonal with respect to the weight function $w_{\alpha,c}$.
At this point, we denote by
\begin{equation} \label{AJ}
A_{n+1}^J(f) = \sum_{i=1}^{n+1} \tilde{w}_i^{(n+1)} f\left(\tilde{x}_i^{(n+1)} \right)    
\end{equation}
and
\begin{equation} \label{AL}
A_{n+1}^L(f) = \sum_{i=1}^{n+1} \tilde{\lambda}_i^{(n+1)} f\left(\tilde{\xi}_i^{(n+1)} \right)    
\end{equation}
the anti-Gaussian quadrature rules (cf. (\ref{A})) relative to $I_n^J$ and $I_n^L$ (see (\ref{GJ}) and (\ref{GL})), respectively.
Then, the corresponding averaged Gauss rules are
\begin{equation*}
    \tilde{A}^J_{2n+1}(f)= \frac{1}{2} \left(I^J_n(f)+A^J_{n+1}(f)\right)
\end{equation*}
and
\begin{equation*}
    \tilde{A}^L_{2n+1}(f)= \frac{1}{2} \left(I^L_n(f)+A^L_{n+1}(f)\right).
\end{equation*}
Finally, the error of the Gaussian quadrature formula (see (\ref{errore totale})) is estimated as 
\begin{equation} \label{err_antiG}
E_{n, \nu,\alpha,c}(f) \approx  \tilde{E}_{n, \nu,\alpha,c}(f):=\tilde{E}_{n, \nu,\alpha,c}^J(f)-\tilde{E}^L_{n,\alpha,c}(f),
\end{equation}
where
\begin{equation*} 
\tilde{E}_{n, \nu,\alpha,c}^J(f) = \tilde{A}_{2n+1}^J(f)-  I_n^J(f) \quad {\rm and} \quad  \tilde{E}_{n,\alpha,c}^L(f) = \tilde{A}_{2n+1}^L(f)-  I_n^L(f).
\end{equation*}

For the generalized Gauss-Laguerre rule the recurrence coefficients and the values of the orthogonal polynomials are explicitly known, and condition (\ref{internal}) holds (see \cite[Theorem 4]{Laurie}). 
Moreover, by using relations (\ref{ab}) and (\ref{polinomi Lag c}) it is not difficult to prove the same result also for the Gaussian rule $I_n^L$.
Therefore, we have that the anti-Gaussian formula $A_{n+1}^L$ is internal, i.e., $\tilde{\xi}_1^{(n+1)} \in [0,+\infty)$.

For what concerns the anti-Gaussian rule $A_{n+1}^J$, we do not have at disposal  the recurrence coefficients and the expressions of the corresponding orthogonal polynomials (see \cite{GaussBessel}). Hence, condition (\ref{internal}) can only be verified numerically. 
In our numerical experiments, independently of $\nu$, $c$ and working with $n=100$, we have observed negative values of $\tilde{x}_1^{(n+1)}$ for $-1 < \alpha< \tilde{\alpha}$, where $\tilde{\alpha} \in (-0.8,-0.7)$.

\section{Construction of generalized averaged Gauss rules} \label{section3}

In this section we describe the generalized averaged Gauss rule $\hat{A}_{2n+1}$, introduced in \cite{S}, associated with a generic Gauss formula $I_n$.
It is a $(2n+1)$-point formula which can be represented by a single symmetric tridiagonal matrix $\hat{J}_{2n+1} \in \mathbb{R}^{(2n+1) \times (2n+1)}$, that is,
\begin{equation*} 
\hat{J}_{2n+1}=
\begin{bmatrix} 
J_n & \sqrt{\beta_n} e_n & \mathbf{0} \\
\sqrt{\beta_n} e_n^T & \alpha_n & \sqrt{\beta_{n+1}} e_1^T \\
\mathbf{0} & \sqrt{\beta_{n+1}} e_1 & J'_n
\end{bmatrix}  ,
\end{equation*}
where $e_1=(1,0,\ldots,0)^T \in \mathbb{R}^n$, $J_n$ is as in (\ref{matriceJ_n}) and $J'_n$ is obtained by reversing the order of the rows and column of $J_n$, that is, 
\begin{equation*} 
J'_n=    \begin{bmatrix}
\alpha_{n-1} & \sqrt{\beta_{n-1}} & & & 0\\
\sqrt{\beta_{n-1}} & \alpha_{n-2} & \sqrt{\beta_{n-2}} \\
 &  \ddots & \ddots & \ddots\\
 & & \sqrt{\beta_2} & \alpha_1 & \sqrt{\beta_1} \\
0 & & & \sqrt{\beta_1} & \alpha_0
\end{bmatrix}  \in \mathbb{R}^{n \times n}.
\end{equation*}
The  generalized averaged Gauss formula can also be described in a more compact form. 
Indeed,  $\hat{A}_{2n+1}$ can be written as
\begin{equation} \label{Ahat}
\hat{A}_{2n+1}(f) = \frac{\beta_{n+1}}{\beta_n+\beta_{n+1}} I_n(f) + \frac{\beta_{n}}{\beta_n+\beta_{n+1}} \Bar{A}_{n+1}(f),
\end{equation}
where the quadrature formula 
\begin{equation} \label{22bis}
    \Bar{A}_{n+1}(f) = \sum_{i=1}^{n+1} \Bar{\sigma}_i^{(n+1)} f \left( \Bar{\tau}_i^{(n+1)} \right)
\end{equation}
is determined by the symmetric tridiagonal matrix $\Bar{J}_{n+1}\in \mathbb{R}^{(n+1) \times (n+1)}$, defined as
\begin{equation}\label{Jbar}
 \Bar{J}_{n+1}=
\begin{bmatrix} 
J_n & e_n \sqrt{\beta_n+\beta_{n+1}} \\
e_n^T \sqrt{\beta_n+\beta_{n+1}} & \alpha_n  \\
\end{bmatrix}, 
\end{equation}
(see \cite{RS21}).
By construction, the $(n+1)$-point quadrature rule (\ref{22bis}) has essentially the same properties of the anti-Gaussian rule $A_{n+1}$, but in this case the nodes are internal, i.e.,  $\Bar{\tau}_1^{(n+1)} \in [0, + \infty)$, if and only if
\begin{equation} \label{condizione2}
    \frac{\pi_{n+1}(0)}{\pi_{n-1}(0)} \geq \beta_{n+1}, \quad n \geq 1,
\end{equation}
(see \cite{S}).

At this point, we denote the generalized averaged Gaussian rules corresponding to $I_n^J$ and $I_n^L$ by
\begin{equation*}
\hat{A}^J_{2n+1}(f) = \frac{\beta^J_{n+1}}{\beta^J_n+\beta^J_{n+1}} I_n^J(f) + \frac{\beta^J_{n}}{\beta^J_n+\beta^J_{n+1}} \Bar{A}^J_{n+1}(f)
\end{equation*}
and
\begin{equation*}
\hat{A}^L_{2n+1}(f) = \frac{\beta^L_{n+1}}{\beta^L_n+\beta^L_{n+1}} I_n^L(f) + \frac{\beta^L_{n}}{\beta^L_n+\beta^L_{n+1}} \Bar{A}^L_{n+1}(f),
\end{equation*}
where the formulas
\begin{equation*}
    \Bar{A}^J_{n+1}(f) = \sum_{i=1}^{n+1} \Bar{w}_i^{(n+1)} f \left( \Bar{x}_i^{(n+1)} \right) \quad {\rm and} \quad \Bar{A}^L_{n+1}(f) = \frac{1}{c^{\alpha+1}} \sum_{i=1}^{n+1} \Bar{\lambda}_i^{(n+1)} f \left( \Bar{\xi}_i^{(n+1)} \right)
\end{equation*}
are associated with the tridiagonal matrices $\Bar{J}_{n+1}^{J}$ and $\Bar{J}_{n+1}^{L}$, obtained by considering in (\ref{Jbar}) the recursion coefficients $\alpha_k^J$, $\beta_k^J$ and $\alpha_k^L$, $\beta_k^L$, respectively.
Finally, the error of the Gaussian quadrature (see (\ref{errore totale})) is estimated as 
\begin{equation} \label{err_averG}
E_{n, \nu,\alpha,c}(f) \approx  \hat{E}_{n, \nu,\alpha,c}(f):=\hat{E}_{n, \nu,\alpha,c}^J(f)-\hat{E}^L_{n,\alpha,c}(f),
\end{equation}
where
\begin{equation*} 
\hat{E}_{n, \nu,\alpha,c}^J(f) = \hat{A}_{2n+1}^J(f)-  I_n^J(f) \quad {\rm and} \quad  \hat{E}_{n, \alpha,c}^L(f) = \hat{A}_{2n+1}^L(f)-  I_n^L(f).
\end{equation*}
From (\ref{condizione2}) and by using again relations (\ref{ab}) and (\ref{polinomi Lag c}), it is trivial to prove that for the generalized averaged Gauss rule $\Bar{A}^L_{n+1}$ the condition for the internality is $\alpha \geq 1$. This means that for $-1<\alpha <1$ the smallest node of $\Bar{A}^L_{n+1}$ is negative.
As before, the behavior of the rule $\Bar{A}_{n+1}^J$ can only be verified numerically. In particular, independently of $\nu$, $c$ and working with $n=100$, we have found negative values of $\Bar{x}_1^{(n+1)}$ for $-1 < \alpha < \Bar{\alpha}$, with $\Bar{\alpha} \in (1,1.1)$.

\section{A tentative a priori estimate} \label{section4}

In the previous sections we have described how, having at disposal the recurrence coefficients of the corresponding orthogonal polynomials, the averaged and generalized averaged Gaussian rules can be easily constructed and employed to approximate the quadrature error $E_{n,\nu,\alpha,c}(f)$ (see (\ref{errore totale})). In this section, by using the results of Barrett \cite{B} regarding the derivation of an asymptotic expression for the error of the Gauss-Laguerre formula, we present a tentative a priori estimate of the error $E_{n,\nu,\alpha,c}(f)$.
In particular, as remarked in the introduction, since $\left\vert J_{\nu} (x)\right\vert \leq 1$, for $\nu \geq 0$, $x \in \mathbb{R}$, we can interpret the weight function $w_{\nu,\alpha,c}$ as a perturbation of the weight of the Gauss-Laguerre rule. 
Therefore, the idea is to employ the error estimate for the Gauss-Laguerre formula to approximate not only $E_{n, \alpha,c}^L(f)$, but also $E_{n,\nu,\alpha,c}^J(f)$ (see (\ref{11bis})).
First, let consider the error of the Laguerre rule for the classical case of $c=1$, that is, $E^L_{n,\alpha,1}(f)$.
For any given $R>1$, the set $\Gamma_R = \left\lbrace z \in \mathbb{C} \mid \Re \sqrt{-z}=\ln R \right\rbrace $ represents a parabola of the complex plane positively oriented, symmetric with respect to the real axis and with vertex at $-\ln R$.
Barrett \cite{B} showed that, if $f(z)$ has no singularities on or within $\Gamma_R$, except for a pair of simple poles $z_0$ and its conjugate $\Bar{z_0}$, then, for $n \rightarrow  \infty$,
\begin{equation} \label{errore Barrett}
    E_{n,\alpha,1}^L(f) \sim  4 \pi \Re \left\lbrace {\rm Res}(f(z),z_0) e^{-i \alpha \pi} \left[ \exp \sqrt{-z_0} \right]^{-2 \sqrt{\Bar{n}}} \right\rbrace,
\end{equation}
where $\Bar{n}=4n+\alpha+2$ and the symbol ${\rm Res(\cdot,\cdot)}$ denotes the residue.
Formula (\ref{errore Barrett}) follows from the fact that $E_{n,\alpha.1}^L(f)$ can be represented as the contour integral
\begin{equation*}
    E_{n,\alpha,1}^L(f) = \frac{1}{2 \pi i} \int_{\Gamma} \frac{q_n^L(z)}{L_n^{(\alpha)}(z)} f(z) dz,
\end{equation*}
where $L_n^{(\alpha)}(z)$ is the generalized Laguerre polynomial, $q_n^L(z)$ is the associated function defined by
\begin{equation*}
    q_n^L(z) = \int_0^{+\infty} \frac{w_{\alpha}(x) L_n^{(\alpha)}(x)}{z-x} dx,
\end{equation*}
with $w_{\alpha}(x)=x^{\alpha} e^{-x}$ and $\Gamma$ is a contour containing $[0, +\infty)$ and such that no singularity of $f(z)$ lies on or within  the contour (see \cite{Szego}).
Then, by using the relation (see  \cite{El})
\begin{equation}
    \frac{q_n^L(z)}{L_n^{(\alpha)}(z)} \sim -2e^{-i\pi \alpha}z^{\alpha}e^{-z}\frac{K_{\alpha}\left(\sqrt{\Bar{n} z e^{-i\pi}}\right)}{I_{\alpha}\left(\sqrt{\Bar{n} z e^{-i\pi}}\right)},
\end{equation}
where $I_{\alpha}$, $K_{\alpha}$ are the modified Bessel functions of order $\alpha$ of the first and second kind, respectively, and the asymptotic formulas (see \cite[p.377, 9.7.1-9.7.2]{Abramowitz})
\begin{equation*}
    I_{\alpha}(t) \sim \frac{e^t}{\sqrt{2 \pi t}}, \;
    K_{\alpha}(t) \sim e^{-t}\sqrt{\frac{\pi}{2t}},
\end{equation*}
valid for large $|t|$, $|\arg (t) |<\pi/2$, we obtain
\begin{equation} \label{asymptotic}
    \frac{q_n^L(z)}{L_n^{(\alpha)}(z)} \sim -2 \pi e^{-i \alpha \pi} z^{\alpha} e^{-z} \left[ \exp \left(\sqrt{-z} \right) \right]^{-2\sqrt{\Bar{n}}}=: \Phi_n(z), \quad z \notin [0, +\infty).
\end{equation}
Finally, formula (\ref{errore Barrett}) can be derived by choosing $\Gamma=\Gamma_R \cup C_1 \cup C_2$, where $C_1$, $C_2$ are two arbitrary small circles surrounding the two poles, and by using relation (\ref{asymptotic}).
After simple computations, it can be verified that, by replacing $w_{\alpha}(x)$ with $w_{\alpha,c}(x)$, one obtains
\begin{equation} \label{errore_parte_reale}
    E_{n,\alpha,c}^L(f) \sim  \frac{4 \pi}{c^{\alpha-1}} \Re \left\lbrace {\rm Res}(f(z),z_0) e^{-i \alpha \pi} \left[ \exp \sqrt{-cz_0} \right]^{-2 \sqrt{\Bar{n}}} \right\rbrace .
\end{equation}

\begin{figure}
\begin{center}
\includegraphics[scale=0.35]{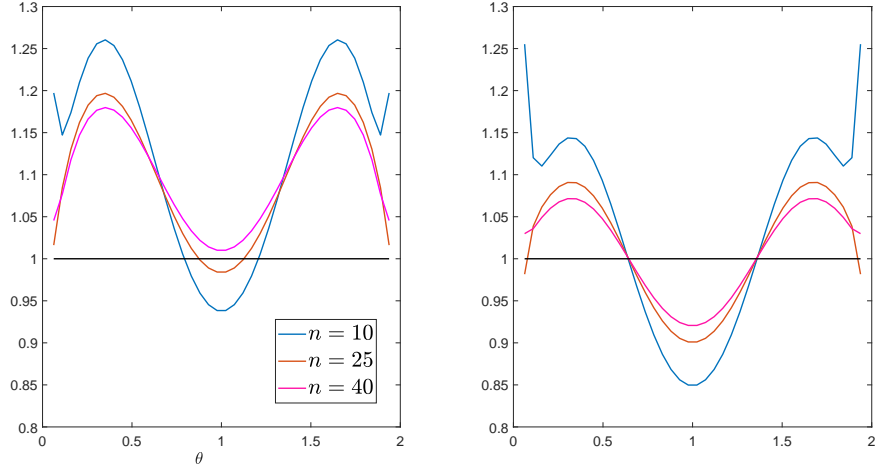}
\end{center}
\caption{The ratios $\frac{\Psi_n^J(z)}{\Phi_n(z)}$ (left) and $\frac{\Psi_n^L(z)}{\Phi_n(z)}$ (right) for $z=r e^{i \theta \pi}$, with $r=4$, $\theta \in ( 0, 2)$. In this case $\nu=1, c=0.5, \alpha=0.3$.} \label{q1}
\end{figure}

\begin{figure}
\begin{center}
\includegraphics[scale=0.35]{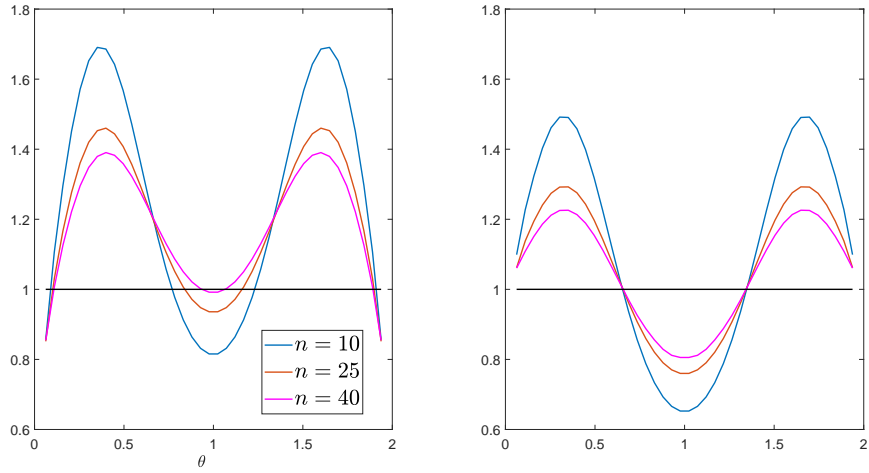}
\end{center}
\caption{The ratios $\frac{\Psi_n^J(z)}{\Phi_n(z)}$ (left) and $\frac{\Psi_n^L(z)}{\Phi_n(z)}$ (right) for $z=r e^{i \theta \pi}$, with $r=4$, $\theta \in ( 0, 2)$. In this case $\nu=0, c=1, \alpha=-0.5$.} \label{q2}
\end{figure}

Now, let consider the error $E_{n,\nu,\alpha,c}^J(f)$. It can also be written as the contour integral 
\begin{equation*}
    E_{n,\nu,\alpha,c}^J(f) = \frac{1}{2 \pi i} \int_{\Gamma} \frac{q_n^J(z)}{p_n(z)} f(z) dz,
\end{equation*}
where $p_n(z)$ is the orthogonal polynomial relative to the weight function $w_{\nu,\alpha,c}$ and
\begin{equation*}
    q_n^J(z) = \int_0^{+\infty} \frac{w_{\nu, \alpha,c}(x) p_n(x)}{z-x} dx.
\end{equation*}
However, as remarked before, we do not have at disposal the analytic expression of $p_n(z)$ and, therefore, an asymptotic formula analog to (\ref{asymptotic}) can not be derived.
In order to justify the use of estimate (\ref{errore_parte_reale}) also for $E_{n,\nu,\alpha,c}^J(f)$, we numerically evaluate the functions
\begin{equation} \label{ratio}
    \Psi_n^J(z) := \frac{q_n^J(z)}{p_n(z)} \quad {\rm and} \quad \Psi_n^L(z) := \frac{q_n^L(z)}{L_n^{(\alpha)}(z)}
\end{equation}
and check if the approximation $\Phi_n(z)$ (cf. (\ref{asymptotic})) can also be used for $\Psi_n^J(z)$.
In Figures \ref{q1}-\ref{q2}, for different values of $n, \nu, \alpha,c$,  we plot the ratios 
\begin{equation*}
    \frac{\Psi_n^J(z)}{\Phi_n(z)} \quad {\rm and} \quad \frac{\Psi_n^L(z)}{\Phi_n(z)},
\end{equation*}
in which $z=r e^{i \theta \pi}$, with $r=4$, $\theta \in ( 0, 2)$.
We remark that $\frac{\Psi_n^L(z)}{\Phi_n(z)} \sim 1$, for $n \rightarrow \infty$, by (\ref{asymptotic}).
The results show that approximation (\ref{asymptotic}) works rather good also for $\Psi_n^J(z)$, and the situation is analog for other values of the parameters.
Hence, the idea is to use the approximation (\ref{errore_parte_reale}) also for $E_{n,\nu,\alpha,c}^J(f)$.
Finally, since $E_{n,\nu,\alpha,c}(f) = E_{n,\nu,\alpha,c}^J(f)-E_{n,\alpha,c}^L(f)$, we consider the estimate
\begin{equation} \label{stimaB}
    \left\vert E_{n,\nu,\alpha,c}(f) \right\vert  \approx 2  \mathcal{E}_{n,\alpha,c}(f), 
\end{equation}
with 
\begin{equation} \label{34bis} 
    \mathcal{E}_{n,\alpha,c}(f) :=  \frac{4 \pi}{c^{\alpha-1}} \left\vert {\rm Res}(f(z),z_0) \left[ \exp \sqrt{-cz_0} \right]^{-2 \sqrt{\Bar{n}}} \right\vert,
\end{equation}
(cf. (\ref{errore_parte_reale})).

\section{Numerical examples} \label{section5}

In this section we present some numerical experiments which confirm the reliability of the error estimates (\ref{err_antiG}), (\ref{err_averG}) and (\ref{stimaB}). We remark that all the computations reported in this work are carried out in Matlab by using high-precision arithmetic, specifically with $120$ significant decimal digits.
The main reason for this choice is that we do not known explicitly the recurrence coefficients $\alpha_k^J$, $\beta_k^J$ and hence we are forced to employ a numerical scheme to derive them. However, it is well known (see e.g. \cite{G}) that this computation can be inaccurate for growing $k$, because the problem is severely ill conditioned. Therefore, even if in \cite{GaussBessel} an alternative more stable approach is presented, the use of high-precision arithmetic allows to considerably increase the number of quadrature points and to obtain an absolute error in the approximations of the order of the machine precision. 
The Matlab routine for the computation of the recurrence coefficients $\alpha_k^J, \beta_k^J$ is taken from  \cite{G2006}, while the code for the implementation of the Gauss-Laguerre quadrature rule from \cite{T}.

\begin{example}
Consider the integral 
\begin{equation*}
    I_{\nu,\alpha,c}(f)=\int_0^{\infty} f(x) x^{\alpha}e^{-cx}J_{\nu}(x) dx,
\end{equation*}
with 
\begin{equation*}
    f(x) = \frac{1}{1+e^{-x}}.
\end{equation*}
In order to use error estimate (\ref{stimaB}), we start by studing the poles of $f$.
A simple analysis shows that they are given by the set
\begin{equation*}
    z_k = -i (\pi +2k \pi), \quad k \in \mathbb{Z},
\end{equation*}
and the closest to the real axis are $z_0$ and its conjugate $z_{1}$, that is, $\pm i \pi$. 
As for the residue (cf, (\ref{errore_parte_reale})), we obtain
\begin{equation*}
    {\rm Res} \left(f(z),z_0\right) = 1.
\end{equation*}
Therefore, by using (\ref{stimaB}) and (\ref{34bis}), we have that
\begin{equation}  \label{apriori1}
     |E_{n,\nu,\alpha,c}(f)| \approx 8 \pi c^{1-\alpha}e^{-\sqrt{2c \pi \Bar{n}}}.
\end{equation}
In Figure \ref{fig1es1} we compare, for different values of $\nu,\alpha,c$, the quadrature error $| E_{n,\nu,\alpha,c}(f) |$, obtained by considering a reference solution, with the approximations $| \tilde{E}_{n,\nu,\alpha,c}(f) |$, $| \hat{E}_{n,\nu,\alpha,c}(f)|$ (cf. (\ref{err_antiG})-(\ref{err_averG})) and estimate (\ref{apriori1}).
We can see the very good agreement between the error and the approximations given by the averaged and generalized averaged Gaussian rules.
Moreover, the examples reveal that the a priori estimate (\ref{apriori1}) is rather accurate.

\begin{figure}
\begin{center}
\includegraphics[scale=0.35]{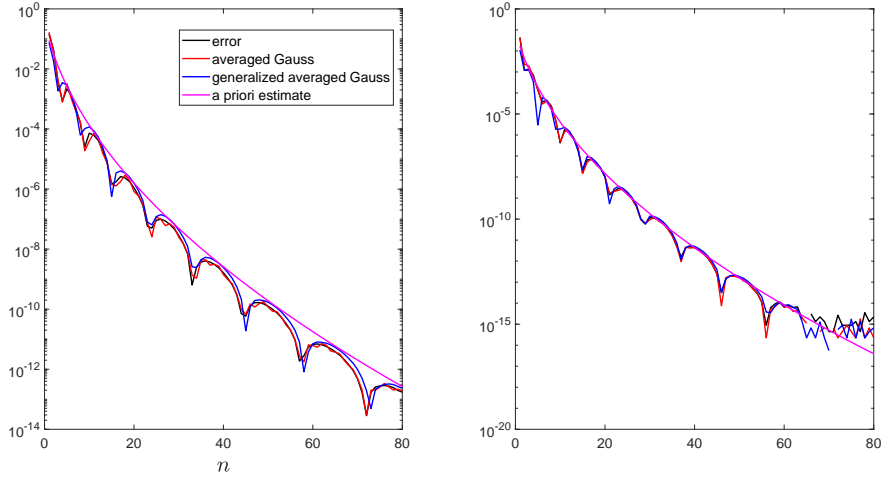}
\end{center}
\caption{The quadrature error and its approximations $| \tilde{E}_{n,\nu,\alpha,c}(f) |$, $| \hat{E}_{n,\nu,\alpha,c}(f)|$ and  (\ref{apriori1}) for $\nu=1, c=0.5, \alpha=1.7$ (left) and $\nu=0.5, c=0.8, \alpha=1.5$ (right).} \label{fig1es1}
\end{figure}
\end{example}

\begin{example}
Consider now the case of
\begin{equation*}
    f(x) = \frac{1}{1+x^2}.
\end{equation*}
This function has only two poles $\pm i$. 
As for the residue (cf. (\ref{errore_parte_reale})), we obtain
\begin{equation*}
    {\rm Res} \left(f(z),i\right) = -\frac{i}{2}.
\end{equation*}
Therefore, by using (\ref{stimaB}) and (\ref{34bis}), we have that
\begin{equation} \label{apriori}
   |E_{n,\nu,\alpha,c}(f)| \approx 4 \pi c^{1-\alpha}e^{-\sqrt{2c \pi \Bar{n}}}.
\end{equation}  
In Figure \ref{fig1es2} we compare, for different values of $\nu,\alpha,c$, the quadrature error $| E_{n, \nu,\alpha,c}(f) |$, obtained by considering a reference solution, with the approximations $| \tilde{E}_{n, \nu,\alpha,c}(f) |$, $| \hat{E}_{n, \nu,\alpha,c}(f) |$ (cf. (\ref{err_antiG})-(\ref{err_averG})) and  estimate (\ref{apriori}).

\begin{figure}
\begin{center}
\includegraphics[scale=0.35]{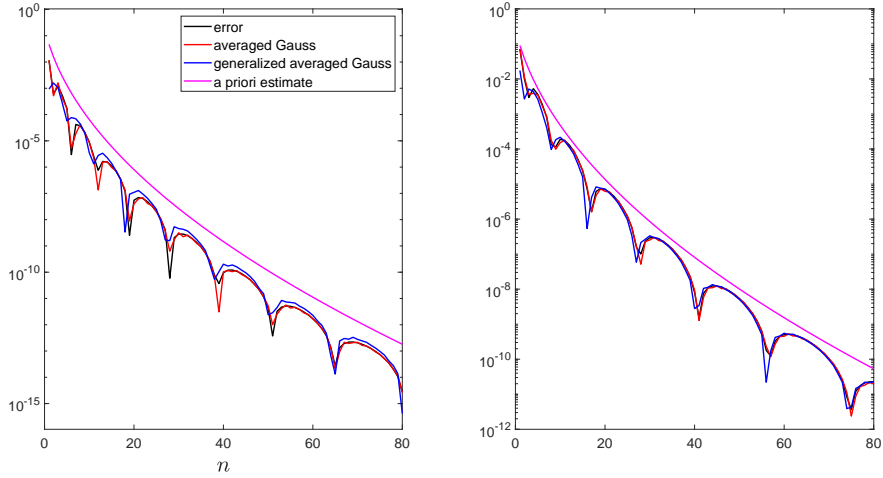}
\end{center}
\caption{The quadrature error and its approximations $| \tilde{E}_{n,\nu,\alpha,c}(f) |$, $| \hat{E}_{n,\nu,\alpha,c}(f) |$ and  (\ref{apriori}) for $\nu=1, c=1.5, \alpha=1$ (left) and $\nu=0, c=1, \alpha=1.5$ (right).} \label{fig1es2}
\end{figure}
\end{example}

\section{Conclusions} \label{section6}

\begin{figure}
\begin{center}
\includegraphics[scale=0.35]{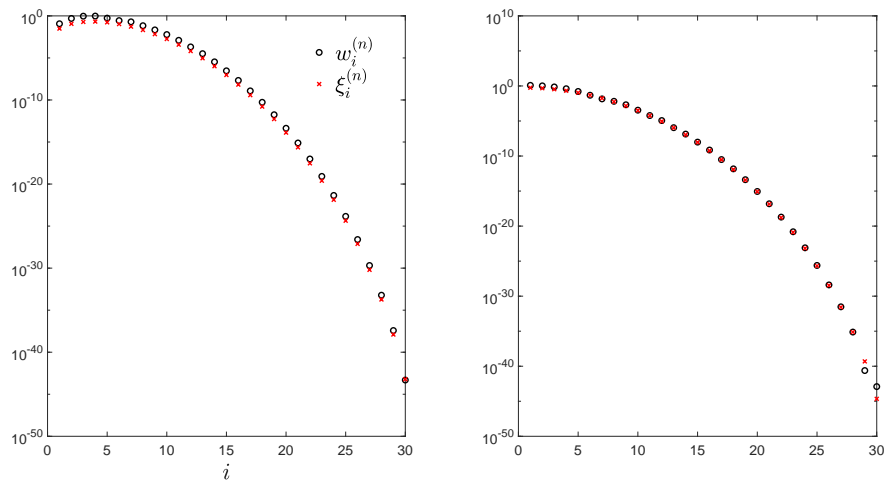}
\end{center}
\caption{The weights $w_i^{(n)}$, $\xi_i^{(n)}$, $i=1,\ldots,n$ (cf. (\ref{GJ})-(\ref{GL})) for $\nu=1$, $\alpha=0.7$, $c=0.5$ (left) and $\nu=0$, $\alpha=-0.5$, $c=0.8$ (right). In both cases $n=30$.} \label{pesi_decay}
\end{figure}

In this work the error estimates of the Gaussian quadrature formula introduced in \cite{GaussBessel} are considered. 
In particular, a posteriori error approximations given by the averaged and generalized averaged Gaussian rules have been constructed and showed to be very accurate.
Moreover, starting from numerical experiments regarding the ratio $q_n^J/p_n$ (see (\ref{ratio})), an heuristic but quite effective a priori error estimate has been introduced.
We remark that, having at disposal an a-priori estimate of the error and by noting that, similar to the Gauss-Laguerre rule, the weights $w_i^{(n)}$, $i=1,\ldots,n$, (cf. (\ref{GJ})) decay exponentially (see Figure \ref{pesi_decay}), a truncated approach can also be introduced, to reduce the number of function evaluation (see \cite{AN}).

\section*{Acknowledgements}

This work was partially supported by GNCS-INdAM and CINECA under HPC-TRES program award number 2019-04. The author is member of the INdAM research group GNCS.

\end{document}